\documentclass[reqno,11pt]{article}

\usepackage{amsmath,latexsym,amssymb,times,mathptm}

\def\sqr#1#2{{\vcenter{\hrule height.#2pt
        \hbox{\vrule width.#2pt height#1pt \kern#1pt
                \vrule width.#2pt}
        \hrule height.#2pt}}}
\def\square{\mathchoice\sqr64\sqr64\sqr{4}3\sqr{3}3}
\def\QED{\hfill$\square$\break}
\def\demo{\noindent{\bf Proof. }}
\def\Ext{{\rm Ext}}
\def\Hom{{\rm Hom}}
\def\ann{{\rm Ann}}
\def\rar{\rightarrow}
\def\lrar{\longrightarrow}
\def\tratto{\mbox{\rule{2mm}{.2mm}$\;\!$}}
\def\Tor{\operatorname{Tor}}
\def\Ass{\operatorname{Ass}}

\def\image{\operatorname{image}}
\def\ann{\operatorname{ann}}
\def\soc{\operatorname{socle}}
\def\o{\overline}

\newtheorem{Theorem}{\bf Theorem}[section]
\newtheorem{Lemma}[Theorem]{\bf Lemma}
\newtheorem{Corollary}[Theorem]{\bf Corollary}
\newtheorem{Proposition}[Theorem]{\bf Proposition}

\newtheorem{Remark}[Theorem]{\bf Remark}
\newtheorem{Example}[Theorem]{\bf Example}

\begin{document}

\baselineskip=12.5pt

\pagestyle{empty}

\ \vspace{1.7in}

\noindent {\LARGE\bf Integral closure of ideals and \\
annihilators of homology}

\vspace{.25in}

\noindent {\large\sc Alberto \ Corso} \\  Department of
Mathematics, University of Kentucky, Lexington, KY 40506, USA
\\
{\it E-mail}: {\tt corso@ms.uky.edu}

\bigskip

\noindent {\large\sc Craig \ Huneke} \\ Department of Mathematics,
University of Kansas, Lawrence, KS 66045, USA \\
{\it E-mail}: {\tt huneke@math.ukans.edu}

\bigskip

\noindent {\large\sc Daniel \ Katz} \\ Department of Mathematics,
University of Kansas, Lawrence, KS 66045, USA \\
{\it E-mail}: {\tt dlk@math.ukans.edu}

\bigskip

\noindent {\large\sc Wolmer \ V. \ Vasconcelos} \\ Department of
Mathematics, Rutgers University, Piscataway, NJ 08854, USA \\
{\it E-mail}: {\tt vasconce@math.rutgers.edu}

\vspace{2.4cm}

\section{Introduction}

Let $(R, {\mathfrak m})$ be a local Noetherian ring. Given an
$R$-ideal $I$ of height $g$, a closely related object to $I$ is
its {\it integral closure} $\overline{I}$. This is the set
$($ideal, to be precise$)$ of all elements in $R$ that satisfy an
equation of the form
\[
X^m+b_1 X^{m-1} + b_2 X^{m-2} + \cdots + b_{m-1} X + b_m = 0,
\]
with $b_j \in I^j$ and $m$ a non-negative integer. Clearly one has
that $I \subseteq \overline{I} \subseteq \sqrt{I}$, where
$\sqrt{I}$ is the {\it radical} of $I$ and consists instead of the
elements of $R$ that satisfy an equation of the form $X^q-b=0$ for
some $b \in I$ and $q$ a non-negative integer. While \cite{EHV}
already provides direct methods for the computation of $\sqrt{I}$,
the nature of $\overline{I}$ is complex. Even the issue of
validating the equality $I = \overline{I}$ is quite hard and
relatively few methods are known \cite{CHV}. In general, computing
the integral closure of an ideal is a fundamental problem in
commutative algebra. Although it is theoretically possible to
compute integral closures, practical computations at present
remain largely out-of-reach, except for some special ideals, such
as monomial ideals in polynomial rings over a field. One known
computational approach is through the theory of Rees algebras: It
requires the computation of the integral closure of the Rees
algebra ${\mathcal R}$ of $I$ in $R[t]$. However, this method is
potentially wasteful since the integral closure of all the powers
of $I$ are being computed at the same time.  On the other hand,
this method has the advantage that for the integral closure
$\overline{A}$ of an affine algebra $A$ there are readily
available {\em conductors}: given $A$ in terms of generators and
relations $($at least in characteristic zero$)$ the Jacobian ideal
${\rm Jac}$ of $A$ has the property that ${\rm Jac} \cdot
\overline{A} \subseteq A$, in other words, $\overline{A}\subseteq
A: {\rm Jac}$. This fact is the cornerstone of most current
algorithms to build $\overline{A}$ \cite{deJ,Vas}.

On a seemingly unrelated level, let $H_i=H_i(I)$ denote the
homology modules of the {\it Koszul complex} ${\mathbb K}_{*}$
built on a minimal generating set $a_1, \ldots, a_n$ of $I$. It is
well known that all the non-zero Koszul homology modules $H_i$ are
annihilated by $I$, but in general their annihilators tend to be
larger. To be precise, this article outgrew from an effort to
understand our basic question:
\begin{quote}
\emph{Are the annihilators of the non-zero Koszul homology modules
$H_i$ of an unmixed ideal $I$ contained in the integral closure
$\overline{I}$ of $I$?}
\end{quote}
\noindent We are particularly interested in the two most
meaningful Koszul homology modules, namely $H_1$ and $H_{n-g}$ ---
the last non-vanishing Koszul homology module. Of course the case
that matters most in dealing with the annihilator of the latter
module is when $R$ is not Gorenstein. We also stress the necessity
of the unmixedness requirement on $I$ in our question. Indeed, let
$R=k[\![x,y,z,w]\!]$ with $k$ a field characteristic zero. The
ideal $I=(x^2-xy,-xy+y^2,z^2-zw,-zw+w^2)$ is an height two mixed
ideal with ${\rm Ann}(H_1) = \overline{I} = (I, xz-yz-xw+yw)$ and
${\rm Ann}(H_2) = \sqrt{I} = (x-y, z-w)$. It is interesting to
note that this ideal has played a significant role in \cite{CHV},
where it was shown that the integral closure of a binomial ideal
is not necessarily binomial, unlike the case of its radical as
shown by Eisenbud and Sturmfels \cite{ES}. A first approach to our
question would be to decide if the annihilators of the Koszul
homology modules are rigid in the sense that the annihilator of
$H_i$ is contained in the annihilator of $H_{i+1}$. Up to radical
this is true by the well-known rigidity of the Koszul complex. If
true, we could concentrate our attention on the last non-vanishing
Koszul homology. Unfortunately, this rigidity is not true. An
example was given by Aberbach: let $R = k[x,y,z]/ (x,y,z)^{n+1}$
and let $E$ be the injective hull of the residue field of $R$.
Then $z$ is in the annihilator of $H_1(x,y;E)$, but $z^n$ does not
annihilate $H_2(x,y;E)$.  It would be good to have an example
where such behavior occurs for the Koszul homology of an ideal on
the ring itself.

An obvious question is: What happens when $I$ is
integrally closed? In Section 2 we provide some validation for
our guiding question. In Corollary~\ref{firstcasecorb}
we show that for any an ${\mathfrak m}$-primary ideal $I$,
that is not a complete intersection, such
that $c \, H_1=0$ and $c \in I \colon {\mathfrak m}$ for $c\in R$,
then $c \in \overline{I}$. In particular, if $I$ is an integrally
closed ideal then ${\rm Ann}(H_1)=I$. We then proceed to study
$\ann(H_1)$ for several classes of ideals with good structure:
these include syzygetic ideals, height two perfect Cohen-Macaulay
ideals, and height three perfect Gorenstein ideals. While in the
case of height two perfect Cohen-Macaulay ideals the Koszul
homology modules are faithful $($see Proposition~\ref{height2}$)$,
in the case of syzygetic ideals we observe that $\ann(H_1)$ can be
interpreted as $I\colon I_1(\varphi)$, where $I_1(\varphi)$ is the
ideal generated by the entries of any matrix $\varphi$ minimally
presenting the ideal $I$ $($see Proposition~\ref{syzygetic}$)$. In
the case of height three perfect Gorenstein ideals we show the
weaker statement that $(\ann(H_1))^2 \subset \overline{I}$ $($see
Theorem~\ref{annGorenstein}$)$.

Section 3 contains variations on a result of Burch, which continue
the theme of this paper in that they deal with annihilators of
homology and integrally closed ideals. The result of Burch that we
have in mind \cite{Burch} asserts that if $\Tor^R_t(R/I,M)$, $M$ a
finitely generated $R$-module, vanishes for two consecutive values
of $t$ less than or equal to the projective dimension of $M$, then
${\mathfrak m}(I\colon{\mathfrak m}) = {\mathfrak m}I$. This has
the intriguing consequence that if $I$ is an integrally closed
ideal with finite projective dimension, then $R_{\mathfrak p}$ is
a regular local ring for all ${\mathfrak p}\in \Ass(R/I)$. In
particular, a local ring is regular if and only if it has an
${\mathfrak m}$-primary integrally closed ideal of finite
projective dimension. A variation of Burch's theorem is given in
Theorem~\ref{strengthen}. We then deduce a number of corollaries.
For instance, we show in Corollary~\ref{cor2} that integrally
closed ${\mathfrak m}$-primary ideals $I$ can be used to test for
finite projective dimension, in the sense that if Tor$_i^R(M,R/I)
= 0$, then the projective dimension of $M$ is at most $i-1$. This
improves Burch's result in that we do not need to assume that two
consecutive Tors vanish. Recent work of Goto and Hayasaka
(\cite{GH1} and \cite{GH2}) has many more results concerning
integrally closed ideals of finite projective dimension.

The annihilator of the conormal module $I/I^2$ is a natural source
of elements in the integral closure of $I$. In Section 4 we study
a class of Cohen-Macaulay ideals whose conormal module is
faithful. We close with a last section giving an equivalent formulation
of our main question, and also include another question which came up in
the course of this study.

\section{Annihilators of Koszul homology}

We start with some easy remarks, that are definitely not sharp
exactly because of their generality. It follows from localization
that $\ann(H_1) \subset \sqrt{I}$. Moreover, for any $R$-ideal $I$
minimally presented by a matrix $\varphi$ we also show that
$\ann(H_1) \subset I \colon I_1(\varphi)$, where $I_1(\varphi)$ is
the ideal generated by the entries of $\varphi$. Things get
sharper when one focuses on the annihilator of the first Koszul
homology modules of classes of ideals with good structural
properties. We conclude the section with a result of Ulrich about
the annihilator of the last non-vanishing Koszul homology module.

\subsection{The first Koszul homology module}

Our first theorem is a general result about annihilators of
Koszul homology. It follows from this theorem that our
basic question has a positive answer for the first Koszul
homology module in the case that $I$ is an integrally
closed $\mathfrak m$-primary ideal.

\begin{Theorem}\label{firstcase}
Let $(R, {\mathfrak m})$ be a local Noetherian ring and let $I$ be
an ${\mathfrak m}$-primary ideal. If $c \in R$ is an element such
that $c \, H_i(I) = 0$ then one of the following conditions hold :
\begin{itemize}
\item[$({\it a})$]
$I : c = \mathfrak m I : c$

\item[$({\it b})$]
There exists $J\subseteq I$ and $x\in R$ such that
$I = J+(cx)$, $\mu(I) = \mu(J)+1$ and $c\ H_i(J) = c\ H_{i-1}(J) = 0$.
\end{itemize}
\end{Theorem}

We will need a lemma before proving Theorem~\ref{firstcase}.

\begin{Lemma}\label{firstcaselem}
Let $J\subseteq R$ be an ideal
and $c, x\in R$. Assume that $(J,cx)$ is primary to the maximal ideal.
Then $\lambda(H_i(J,c)) = \lambda(\ann_cH_i(J,cx))$.
\end{Lemma}
\demo
Induct on $i$. Suppose $i = 0$. The desired equality
of lengths follows immediately from the exact sequence
\[
0 \rar ((J,cx):c)/(J,cx) \lrar R/(J,cx) \stackrel{\cdot c}{\lrar}
R/(J,cx) \lrar R/(J,c) \rar 0.
\]
Suppose $i >0$ and the lemma holds for $i-1$. We have an exact sequence
\[
0 \rar H_i(J,cx)/cH_i(J,cx) \lrar H_i(J,cx,c) \lrar
\ann_c(H_{i-1}(J,cx) \rar 0.
\]
But $H_i(J,cx,c) = H_i(J,c)\oplus H_{i-1}(J,c)$, so
\[
\lambda(H_i(J,c)) + \lambda(H_{i-1}(J,c)) =
 \lambda(\ann_c(H_{i-1}(J,c)) + \lambda(H_i(J,cx)/cH_i(J,cx)).
\]
Using the induction hypothesis, we obtain
$\lambda(H_i(J,c)) = \lambda(H_i(J,cx)/cH_i(J,cx)) =
\lambda(\ann_c H_i(J,cx))$. \QED

\noindent {\bf Proof of Theorem~\ref{firstcase}.} Suppose $({\it
a})$ does not hold. Then there exists $x\in \mathfrak m$ such that
$cx$ is a minimal generator of $I$. We can write $I = J + (cx)$,
for an ideal $J\subseteq I$ satisfying $\mu(I) = \mu(J)+1$. On the
one hand, from the exact sequences
\[
0 \rar H_i(J)/cH_i(J) \lrar H_i(J,c) \lrar \ann_c H_{i-1}(J) \rar
0
\]
and
\[
0 \rar H_i(J)/cxH_i(J) \lrar H_i(J,cx) \lrar \ann_{cx} H_{i-1}(J)
\rar 0
\]
we get
\[
\lambda(H_i(J,c)) = \lambda(H_i(J)/cH_i(J)) + \lambda(\ann_c
H_{i-1}(J))
\]
and
\[
\lambda(H_i(J,cx)) = \lambda(H_i(J)/cxH_i(J)) + \lambda(\ann_{cx}
H_{i-1}(J)).
\]
On the other hand,
\[
\lambda(H_i(J)/cxH_i(J))\geq \lambda(H_i(J)/cH_i(J)) \quad {\rm
and} \quad \lambda (\ann_{cx} H_{i-1}(J)) \geq \lambda(\ann _c
H_i(J)).
\]
Since $cH_i(J,cx) = 0$, $H_i(J,cx) = \ann_c H_i(J,cx)$, so
$\lambda(H_i(J,cx)) = \lambda(H_i(J,c))$, by
Lemma~\ref{firstcaselem}. It follows from this that
$\lambda(H_i(J)/cH_i(J)) = \lambda(H_i(J)/cxH_i(J))$. Thus,
$cH_i(J) = cxH_i(J)$, so $cH_i(J) = 0$, by Nakayama's lemma.
Similarly, since
\[
\lambda(\ann_c H_{i-1}(J)) =
\lambda(\ann_{cx} H_{i-1}(J)),
\]
it follows that $\lambda(H_{i-1}(J)/cH_i(J)) =
\lambda(H_{i-1}(J)/cxH_i(J))$, so $cH_{i-1}(J) = 0$, as before.
\QED

\begin{Corollary}\label{firstcasecora}
Let $(R, {\mathfrak m})$ be a local Noetherian ring and let $I$ be
an ${\mathfrak m}$-primary ideal. If $c\cdot H_1(I) = 0$, then
$I : c = {\mathfrak m}I : c$.
\end{Corollary}
\demo
If $I : c$ properly contains $mI : c$,
then by Theorem~\ref{firstcase}, there exists $J\subseteq I$ and
$x \in \mathfrak m$ such that $I = J + (cx)$, $\mu(I) = \mu(J)+1$
and $c\cdot H_0(J) = 0$. But then, $c \in J$, so $I = J$, a contradiction.
\QED

\begin{Corollary}\label{firstcasecorb}
Let $(R, {\mathfrak m})$ be a local Noetherian ring and let $I$ be
an ${\mathfrak m}$-primary ideal. If $c \in R$ is an element such that
$c\cdot H_1(I) = 0$ and $c \in I : \mathfrak m$, then
$c\in \overline{I}$.
\end{Corollary}
\demo Since $\mathfrak m\subseteq I : c$, we have $\mathfrak m
c\subseteq \mathfrak m I$, by Corollary~\ref{firstcasecora}. By
the determinant trick, $c \in \overline{I}$. \QED

\begin{Corollary}\label{int-closed}
Let $(R, {\mathfrak m})$ be a local Noetherian ring and let $I$ be
an integrally closed ${\mathfrak m}$-primary ideal. Then ${\rm
Ann}(H_1)=I$.
\end{Corollary}
\demo Suppose $\ann H_1(I)$ properly contains $I$. Take $c \in
(\ann H_1(I) \backslash I)\cap (I : \mathfrak m)$. By
Corollary~\ref{firstcasecorb}, $c \in \overline{I} = I$, a
contradiction. Thus, $\ann H_1(I) = I$. \QED

\noindent{\bf Syzygetic ideals:} It follows from the
determinant trick that the annihilator of $I^m/I^{m+1}$ is
contained in $\overline{I}$ for all $m$. Hence, another piece of
evidence in support of our question is given by the close
relationship between $H_1$ and the conormal module $I/I^2$. This
is encoded in the exact sequence
\[
0 \rar \delta(I) \lrar H_1 \lrar (R/I)^n \lrar I/I^2 \rar 0,
\]
where $\delta(I)$ denotes the kernel of the natural surjection
from the second symmetric power ${\rm Sym}_2(I)$ of $I$ onto
$I^2$, ${\rm Sym}_2(I) \twoheadrightarrow I^2$, see \cite{SV}. We
will exploit this exact sequence in at least two occasions:
Proposition~\ref{syzygetic} and Theorem~\ref{faithful}. We recall
that the ideal $I$ is said to be {\it syzygetic} whenever
$\delta(I)=0$.

\begin{Proposition}\label{syzygetic}
Let $R$ be a Noetherian ring. For any $R$-ideal $I$ minimally
presented by a matrix $\varphi$,  $\ann(H_1) \subset I \colon
I_1(\varphi)$, where $I_1(\varphi)$ denotes the ideal generated by
the entries of $\varphi$. If, in addition, $I$ is syzygetic then
$\ann(H_1) = I \colon I_1(\varphi)$.
\end{Proposition}
\demo Let $Z_1$ and $B_1$ denote the modules of cycles and
boundaries respectively. If $x \in \ann(H_1)$ one has that for $z
\in Z_1$ the condition $xz \in B_1$ means that each coordinate of
$z$ is conducted into $I$ by $x$. Thus $x \in I : I_1(\varphi)$.
The reverse containment holds if $I$ is syzygetic. In fact, in
this situation one actually has that $H_1 \hookrightarrow
(I_1(\varphi)/I)^n$. Thus $I \colon I_1(\varphi) \subset {\rm
Ann}(H_1)$. \QED

\begin{Corollary}\label{power-levin}
Let $R$ be a local Noetherian ring, and let $I$ be an ideal of
finite projective dimension $n$. Then $(\ann(H_1))^{n+1}\subseteq
I$.
\end{Corollary}
\demo Assume $I$ is minimally
presented by a matrix $\varphi$. By the above proposition,
$\ann(H_1) \subset I \colon
I_1(\varphi)$. The result then follows immediately from the following
proposition of G. Levin (unpublished). The proof follows from a careful
analysis of Gulliksen's Lemma, 1.3.2 in \cite{GL}.
\QED

\begin{Proposition}\label{levin}
Let $R$ be a local Noetherian ring and let $I$ be an ideal of finite
projective dimension $n$, minimally presented by a matrix $\varphi$.
Then $(I \colon I_1(\varphi))^{n+1}\subseteq I$.
\end{Proposition}

\medskip

\begin{Remark}
{\rm In general, the ideal $I : I_1(\varphi)$ may be larger than
the integral closure of $I$. For example the integrally closed
$R$-ideal $I = (x, y)^2$, where $R$ is the localized polynomial
ring $k[x,y]_{(x,y)}$, is such that $I \colon I_1(\varphi) = (x,
y)$. However, Levin 's proposition shows
that $(I \colon I_1(\varphi))^2 \subset I$. }
\end{Remark}

\medskip

\noindent {\bf Height two perfect ideals:} The first case to
tackle is the one of height two perfect ideals in local
Cohen-Macaulay rings. However the Cohen-Macaulayness of the
$H_i$'s gets into the way. Indeed we have the following fact:

\begin{Proposition}\label{height2}
Let $R$ be a local Cohen-Macaulay ring and let I be a height two
perfect $R$-ideal. Then for all $i$ $($with $H_i \neq 0)$ one has
${\rm Ann}(H_i) = I$.
\end{Proposition}
\demo Consider the resolution of the ideal $I$
\[
0 \rar R^{n-1} \lrar R^n \lrar I \rar 0.
\]
The submodule of $1$-cycles of ${\mathbb K}_{*}$, $Z_1$, is the
submodule $R^{n-1}$ of this resolution. Also, for all $i$ one has
$Z_i = \bigwedge^i Z_1$. All these facts can be traced to
\cite{AH}. This implies that for any $i \leq n-2$, $H_1^{i} = H_i$
--- this multiplication is in $H_{*}({\mathbb K})$. Thus the
annihilator of $H_1$ will also annihilate, say, $H_{n-2}$. But
this is the canonical module of $R/I$, and its annihilator is $I$.
The conclusion now easily follows. \QED

\medskip

\noindent {\bf Gorenstein ideals:} Let us consider a perfect
${\mathfrak m}$-primary Gorenstein ideal in a local Noetherian
ring $R$. In this situation, if $I$ is Gorenstein but not a
complete intersection then ${\rm Ann}(H_1) \neq I$. Otherwise,
$R/I$ would be a submodule of $H_1$. By a theorem of Gulliksen
\cite{GL}, if $H_1$ has a free summand then it must be a complete
intersection. Actually, using Gulliksen's theorem one shows that
if $I$ is ${\mathfrak m}$-primary, Gorenstein but not a complete
intersection, then the socle annihilates $H_1$. Combining
Proposition~\ref{syzygetic} and the work of \cite{CHV} yields the
following result:

\begin{Proposition}
Let $(R, {\mathfrak m})$ be a local Noetherian ring with embedding
dimension at least $2$ and let $I$ be an ${\mathfrak m}$-primary
ideal contained in ${\mathfrak m}^2$ with $R/I$ Gorenstein.
Suppose further that $I$ is minimally presented by a matrix
$\varphi$ and that $I_1(\varphi)={\mathfrak m}$, where
$I_1(\varphi)$ denotes the ideal generated by the entries of
${\varphi}$. Then $\ann(H_1) \subset \overline{I}$.
\end{Proposition}
\demo By Proposition~\ref{syzygetic} and our assumption we obtain
that $\ann(H_1) \subset I \colon I_1(\varphi)=I \colon {\mathfrak
m}$. Our assertion now follows from Lemma~3.6 in \cite{CHV} since
$(I \colon {\mathfrak m})^2=I(I \colon {\mathfrak m})$. \QED

\medskip

For an height three perfect Gorenstein ideal $I$ we have some
evidence that \ $(\ann(H_1))^2 = I \cdot \ann(H_1)$. \ If this
were to hold in general, it would imply that $I \subsetneq
\ann(H_1) \subset \overline{I}$. Thus far, we can prove the weaker
result that the square of the annihilator of $H_1$ is in the
integral closure of $I$.

\begin{Theorem}\label{annGorenstein}
Let $R$ be a local Noetherian ring with ${\rm char}(R) \not=2$ and
let $I$ \/ be a height three perfect Gorenstein ideal minimally
generated by $n \geq 5$ elements. Then
\[
(\ann(H_1))^2 \subset \overline{I}.
\]
\end{Theorem}
\demo Let $a_1, \ldots a_n$ denote a set of minimal generators of
$I$. Notice that $B_1$ and $Z_1$ are submodules of $R^n$ of rank
$n-1$; in general, if $E$ is a submodule of $R^n$ of rank $r$, we
denote by $\det(E)$ the ideal generated by the $r \times r$ minors
of the matrix with any set of generators of $E$ $($as  elements of
$R^n)$.

Let $c \in R$ be such that $cZ_1 \subset B_1$. It suffices to
prove that $c^2 \in \overline{I}$ since the square
of an ideal is always integral over the ideal generated by the squares
of its generators.
Note that $c^{n-1}\det(Z_1) \subset \det(B_1)$. Let $V$ be a valuation
overring of $R$ with valuation $v$; the ideal $IV$ is now
principal and generated by one of the original generators, say
$a_1=a$. By the structure theorem of Buchsbaum and Eisenbud \cite{BE}, we
may assume that $a$ is one of the maximal Pfaffians of the matrix presenting
$I$. Since $I$ is generated by $a$, $B_1V$ is generated by the Koszul
syzygies $(a_2, -a, 0, \ldots, 0), (a_3, 0, -a, \ldots, 0),
\ldots, (a_n, 0, 0, \ldots, -a)$. Hence $\det(B_1V) = (a^{n-1})=
I^{n-1}V$. As for $Z_1V$, one has that $\det(Z_1V)$ includes the
determinant of the minor defining $a^2$ $(a$ is the Pfaffian of
the submatrix$)$. Thus \ $c^{n-1}I^2V \subset I^{n-1}V$, \ which
yields that \ $c^{n-1} \in I^{n-3}V$, as cancellation holds.
Hence, we have that $(n-1)v(c)=v(c^{n-1}) \geq v(I^{n-3}V) =
(n-3)v(IV)$. Finally, this yields
\[
v(c^2) \geq 2\frac{n-3}{n-1} v(IV) \geq v(IV)
\]
and, in conclusion, $c^2 \in \overline{I}$. \QED
\medskip
\begin{Remark}
{\rm It is worth remarking that the above proof shows much more. Recall that
$\overline{I^{\frac{a}{b}}}$ denotes the integral closure of the
ideal generated by all $x\in R$ such that $x^b\in I^a$. From \cite{BE}, we
know that $n = 2k+1$ must be odd. Our proof
shows that $$(\ann(H_1)) \subset \overline{I^{\frac{k-1}{k}}}.$$
As $k$ gets large this is very close to our main objective, proving
that $(\ann(H_1)) \subset \overline{I}$.}
\end{Remark}

\subsection{Last non-vanishing Koszul homology module}

Let us turn our attention towards the tail of the Koszul complex.

\begin{Proposition}
Let $R$ be a one-dimensional domain with finite integral closure.
Then any integrally closed ideal is reflexive. In particular, for
any ideal $I$ its bidual $(I^{-1})^{-1}$ is contained in its
integral closure $\overline{I}$.
\end{Proposition}

\demo We may assume that $R$ is a local ring, of integral closure
$B$.  An ideal $L$ is integrally closed if $L = R \cap LB$. Since
$B$ is a principal ideals domain, $LB =xB$ for some $x$. We claim
that $xB$ is reflexive. Let $C = B^{-1}=\Hom_R(B,R)$ be the
conductor of $B/R$. $C$ is also an ideal of $B$, $C = yB$, and
therefore $C^{-1}= y^{-1}B^{-1}=y^{-1}C$, which shows that
$C^{-1}= B$. This shows that $(L^{-1})^{-1}\subset
(R^{-1})^{-1}\cap ((xB)^{-1})^{-1}= L$. The last assertion follows
immediately by setting $I \subset L = \overline{I}$. \QED

We can interpret the above result as an annihilation of Koszul
cohomology. Let $I= (a_1, \ldots, a_m)$ and let $\mathbb{K}^{*}$
denote the Koszul complex
\[ 0
\rar R \lrar R^m \lrar \bigwedge^2R^m \lrar \cdots \lrar
\bigwedge^m R^m \rar 0,
\]
with differential $\partial (w)= z \wedge w$, where $z = a_1e_1 +
\cdots + a_ne_m$. One sees that $Z^1 = I^{-1}z$, and $B^1= Rz$.
Thus $(I^{-1})^{-1}$ is the annihilator of $H^1$. On the other
hand $H^1 \cong H_{m-1} \cong {\rm Ext}_R^1(R/I, R)$. Let us raise
a related issue: $(I^{-1})^{-1}$ is just the annihilator of
$\Ext_R^1(R/I,R)$, so one might want to consider the following
question which is obviously relevant only if the ring $R$ is not
Gorenstein. Let $R$ be a Cohen-Macaulay geometric integral domain
and let $I$ be a height unmixed ideal of codimension $g$. Is \
$\ann (\Ext_R^g(R/I,R)) = \ann (H_{n-g})$ \ always contained in
$\overline{I}$? Notice that the annihilator of the last non-vanishing
Koszul homology can be identified with $J:(J:I)$ for $J$ an ideal
generated by a maximal regular sequence inside $I$. This follows since the last
non-vanishing Koszul homology is isomorphic to $(J:I)/J$.

\bigskip

We thank Bernd Ulrich for allowing us to reproduce the following
result  \cite{B}, which grew out of conversations
at MSRI $($Berkeley$)$:

\begin{Theorem}[Ulrich]\label{ulrich}
Let $(R,{\mathfrak m})$ be a Cohen-Macaulay local ring, let  $I$
be an ${\mathfrak m}$-primary ideal and let $J\subset I$ be a
complete intersection. Then $J:(J:I) \subset \overline{I}$.
In particular the annihilator of the last non-vanishing Koszul
homology of $I$ is contained in the integral closure of $I$.
\end{Theorem}
\demo We may assume that ${\rm ht}\, J = {\rm ht}\, I$. We may
also assume that $R$ has a canonical module $\omega$. We first
prove:

\begin{Lemma}\label{lemma}
Let $A$ be an Artinian local ring with canonical module
$\omega$ and let $I\subset A$ be an ideal. Then
$0:_{\omega}(0:_AI)= I\omega$.
\end{Lemma}
\demo Note that $0:_{\omega}(0:_AI)= \omega_{R/0:I}$. To show
$I\omega = \omega_{R/0:I}$ note that the socle of $I\omega$ is
1-dimensional as it is  contained in the socle of $\omega$. Hence
we only need to show that $I\omega$ is faithful over $R/0:I$. Let
$x\in \ann_RI\omega$, then  $xI\omega=0$, hence $xI=0$, hence
$x\in 0:I$. \QED

\noindent Returning to the proof of Theorem~\ref{ulrich}, it
suffices to show that $(J:(J:I))\omega \subset I\omega$. But
$(J:(J:I))\omega \subset J\omega:_{\omega}(J:_RI)$. So it suffices
to show $J\omega:_{\omega}(J:_RI) \subset I\omega$. Replacing\,
$R,\omega$\, by\, $R/J$, $\omega_{R/J}= \omega/J\omega$\, we have
to show $0:_{\omega}(0:_{R}I)\subset I\omega$, which holds by
Lemma~\ref{lemma}. \QED

\section{Variations on a theorem of Burch}

Theorem~\ref{strengthen} below is a variation of Burch's theorem
mentioned in the introduction, and strengthens it in the case $I$
is integrally closed. We then deduce a number of corollaries.

\begin{Theorem}\label{strengthen}
Let $(R,{\mathfrak m})$ be a local ring, $I$ an integrally closed
$R$-ideal having height greater than zero and $M$ a finitely
generated $R$-module. For $t\geq 1$, set $J_t :=
\ann(\Tor_t(R/I,M))$. Let $({\mathbb F}_{*}, \varphi_i)$ be a
minimal free resolution of $M$. If $\image(\varphi_t)$ is
contained in $\o{{\mathfrak m}J_t}F_{t-1}$, then
\[
\image(\varphi_t\otimes_R 1_{R/I}) \cap \soc(F_{t-1}/IF_{t-1}) =
0.
\]
\end{Theorem}
\demo Take $u\in F_{t-1}$ such that its residue class modulo $I$
belongs to
\[
\image(\varphi_t\otimes_R 1_{R/I}) \cap \soc(F_{t-1}/IF_{t-1}).
\]
Then $u = \varphi_t(v) + w$, for $v\in F_t$ and $w\in IF_{t-1}$.
For all $x\in {\mathfrak m}$, $\varphi_t(xv) \equiv 0$ modulo
$IF_{t-1}$. Thus for all $j\in J_t$, there exists $z\in F_{t+1}$
such that $\varphi_{t+1}(z) \equiv jxv$ modulo $IF_t$. It follows
that we can write $jxv = \varphi_{t+1}(z) + w_0$, for $w_0\in
IF_t$. Therefore, $jxu = \varphi_t(jxv) + jxw = \varphi_t(w_0) +
jxw$. By hypothesis, we get $jxu \in \o{mJ_t} IF_t$, for all $j\in
J_t$ and all $x\in {\mathfrak m}$. Therefore, by cancellation, $u
\in I_a F_{t-1}$. But since $I$ is integrally closed, $u \in
IF_{t-1}$, which gives what we want. \QED

In the following corollaries, we retain the notation from
Theorem~\ref{strengthen}.

\begin{Corollary}\label{cor1}
Suppose $I$ is integrally closed and ${\mathfrak m}$-primary. If
$\image(\varphi_t)$ is contained in $\o{{\mathfrak m}J_t}F_{t-1}$,
then:
\begin{itemize}
\item[$({\it a})$]
$\image(\varphi_t) \subseteq IF_{t-1}$.

\item[$({\it b})$]
$J_tF_t \subseteq \image(\varphi_{t+1})$.
\end{itemize}
\end{Corollary}
\demo For $({\it a})$, if $\image(\varphi_t\otimes 1_{R/I})$ were
not zero, then it would contain a non-zero socle element, since
$I$ is ${\mathfrak m}$-primary. This contradicts
Theorem~\ref{strengthen}, so $({\it a})$ holds.

\noindent For $({\it b})$, it follows from $({\it a})$ that
$\Tor_t(R/I,M) = F_t/(\image(\varphi_{t+1})+IF_t)$, so $J_tF_t$ is
contained in $\image(\varphi_{t+1}) +IF_t$, and $({\it b})$
follows via Nakayama's Lemma. \QED

The next corollary shows that integrally closed ${\mathfrak
m}$-primary ideals can be used to test for finite projective
dimension.

\begin{Corollary}\label{cor2}
Suppose that $I$ is integrally closed and ${\mathfrak m}$-primary.
Then $M$ has projective dimension less than $t$ if and only if\/
$\Tor_t(R/I,M) = 0$.
\end{Corollary}
\demo The hypothesis implies that $J_t = R$. Therefore,
$\image(\varphi_t)$ is automatically contained in $\o{{\mathfrak
m}J_t}F_{t-1}$. By part $({\it b})$ of Corollary~\ref{cor1}, $F_t
\subseteq \image(\varphi_{t+1})$, so $F_t = 0$, by Nakayama's
Lemma. \QED

\begin{Corollary}
Let $J\subseteq R$ be an ideal and $I$ an integrally closed
${\mathfrak m}$-primary ideal. If $J\subseteq \o{{\mathfrak
m}(IJ:I\cap J)}$, then $J\subseteq I$.
\end{Corollary}
\demo Apply Corollary~\ref{cor1} with $M = R/J$ and $t = 1$. \QED

\begin{Corollary}
Suppose that $R$ is reduced and $M$ has infinite projective dimension
over $R$. Then for all $t\geq 1$, the entries of $\varphi_t$ do
not belong to $\o{{\mathfrak m}\cdot \ann(M)}$. In particular,
each map in the minimal
resolution of $k$ has an entry not belonging to $\o{{\mathfrak
m}^2}$.
\end{Corollary}
\demo Let $I$ be any ${\mathfrak m}$-primary integrally closed
ideal. If the entries of $\varphi_t$ belong to $\o{{\mathfrak
m}\cdot \ann(M)}$, then $\image(\varphi_t)$ is contained in
$\o{{\mathfrak m}J_t}F_{t-1}$. By Corollary~\ref{cor1},
$\image(\varphi_t)$ is contained in $IF_{t-1}$. But the
intersection of the integrally closed ${\mathfrak m}$-primary
ideals is zero, therefore, $\image(\varphi_t) = 0$, contrary to
the hypothesis on $M$. Thus, the conclusion of the corollary
holds.

The last statement follows in the case $R$ is regular from the
fact that the Koszul complex on a minimal set of generators of the
maximal ideal gives a resolution of $k$. If $R$ is not regular, then
$k$ has infinite projective dimension, and the result follows at once
from the first statement. \QED

In regard to the above corollary,
it is well-known that the Koszul
complex of a minimal set of generators of the maximal ideal is part
of a minimal resolution of $k$ in all cases, so for the maps occuring
in the minimal resolution up to  the  dimension of
the ring, the last statement is clear. The new content of the last statement is for the maps past
the dimension of the ring.

\begin{Corollary}
Suppose $I$ is integrally closed and ${\mathfrak m}\in \Ass(R/I)$.
If $\image(\varphi_t)$ is contained in $\o{{\mathfrak
m}J_t}F_{t-1}$, e.g., $\Tor_t(R/I,M) = 0$, then either $M$ has
projective dimension less than $t-1$ or ${\mathfrak m}\in
\Ass(\Tor_{t-1}(R/I,M))$.
\end{Corollary}
\demo Suppose $M$ has projective dimension greater than or equal
to $t-1$. Then $F_{t-1} \not = 0$. By hypothesis, the socle of
$F_{t-1}/IF_{t-1}$ is non-zero, so a non-zero element $u$ in this
socle goes to zero under $\varphi_{t-1}\otimes 1_{R/I}$. But the
theorem implies that the image of $u$ in $\Tor_{t-1}(R/I,M)$
remains non-zero, so the result holds. \QED

\section{The conormal module}

We end the article with a result in the spirit of our
investigation. More precisely we show that the conormal module
$I/I^2$ is faithful for a special class of Cohen-Macaulay ideals.

\begin{Theorem}\label{faithful}
Let $(R, {\mathfrak m})$ be a Gorenstein local ring and $I$ a
Cohen-Macaulay almost complete intersection. Let $\varphi$ be a
matrix minimally presenting $I$. If $I_1(\varphi)$ is a complete
intersection, then $I/I^2$ is a faithful $R/I$-module.
\end{Theorem}
\demo Let $g$ denote the height of $I$, write $n = g+1$ for the
minimal number of generators of $I= (a_1, \ldots, a_n)$. We may assume
 that the ideals generated by any $g$ of the $a_i$'s are
complete intersection ideals. Let $e_i$, with $1 \leq i \leq n$,
denote the $n$-tuple $(0, \ldots, 0,1,0, \ldots, 0)$ where $1$ is
in the $i$-th position. Finally, note that $H_1$ is the canonical
module of $R/I$.

Let us consider the exact sequence
\[
0 \rar \delta(I) \lrar H_1 \stackrel{\theta}{\lrar} (R/I)^n
\stackrel{\pi}{\lrar} I/I^2 \rar 0,
\]
where $\delta(I)$ is the kernel of the natural surjection ${\rm
Sym}_2(I) \twoheadrightarrow I^2$, see \cite{SV}. Notice that for
any $\varepsilon' = \sum r'_j e_j + B_1 \in H_1$, where $\sum r'_j
a_j = 0$, one has $\theta(\varepsilon') = (r'_1 + I)e_1 + \ldots +
(r'_n +I)e_n$ while for any element in $(R/I)^n$ one has $\pi((r_1
+ I)e_1 + \ldots + (r_n + I)e_n)= r_1a_1 + \ldots + r_na_n + I^2$.
Apply $(\tratto)^{\vee}= \Hom_{R/I}(\tratto, H_1)$ to the above
exact sequence. We obtain
\[
0 \rar (I/I^2)^\vee \stackrel{\pi^\vee}{\lrar} \Hom((R/I)^n, H_1)
\stackrel{\theta^\vee}{\lrar} \Hom(H_1,H_1) = R/I \lrar
\delta(I)^\vee \rar 0.
\]
To conclude it will be enough to show that $(I/I^2)^\vee$ is
faithful.

First, we claim that the image of $\theta^\vee$ belongs to
$I_1(\varphi)/I$. In fact, any element of $\Hom((R/I)^n, H_1)$ can
be written as a combination of elementary homomorphism of the form
\[
\xi_i((1+I)e_i)= \varepsilon \qquad \qquad \xi_i((1+I)e_j)=0,
\quad {\rm if} \ i \not= j,
\]
with $\varepsilon = \sum r_j e_j + B_1 \in H_1$, where $\sum r_j
a_j = 0$. Thus, for any $\varepsilon' \in H_1$ we have
\[
(\theta^\vee(\xi_i))(\varepsilon') = \xi_i(\theta(\varepsilon')) =
\xi_i(\sum (r'_j+I) e_j) = (r_i' +I) \varepsilon.
\]
Observe that $(r_i'+I)\varepsilon = (r_i+I)\varepsilon'$ in $H_1$.
Indeed, $r_i'\varepsilon - r_i\varepsilon' = \sum
(r'_ir_j-r_ir'_j)e_j + B_1$. But $\sum (r'_ir_j-r_ir'_j)e_j$ is a
syzygy of the complete intersection $(a_1, \ldots, \widehat{a}_i,
\ldots, a_n)$ and thus it is a Koszul syzygy of the smaller ideal:
hence it is in $B_1$. In conclusion, $\theta^\vee(\xi_i)$ is
nothing but multiplication by $r_i+I \in I_1(\varphi)/I$. Given
that $\varepsilon$ and $i$ were chosen arbitrarily one has that
the image of $\theta^\vee$ is $I_1(\varphi)/I$.

Notice that the number of generators of $I_1(\varphi)$ is strictly
smaller than $n$. So we can say that the image of $\theta^\vee$ is
given say by $(\theta^\vee(\xi_2), \ldots, \theta^\vee(\xi_n))$.
Write, for some $c_i \in R/I$,
\[
\theta^\vee(\xi_1) = \sum_{i\geq 2} c_i \, \theta^\vee(\xi_i).
\]
Hence $\xi_1 - \sum_{i\geq 2} c_i \, \xi_i \in {\rm Ker}
(\theta^\vee) = {\rm Im}(\pi^\vee)$ so that we can find $\gamma
\in (I/I^2)^\vee$ such that
\[
\xi_1 - \sum_{i\geq 2} c_i \, \xi_i = \pi^\vee(\gamma) = \gamma
\circ \pi.
\]
The restriction of these homomorphisms to the first component of
$\Hom ((R/I)^n, H_1)$ gives an homomorphism from $R/I$ to $H_1$.
Now, something that annihilates $\gamma$ would also annihilate the
restriction, but that restriction is faithful. \QED

\medskip

\begin{Remark}
{\rm From the proof of Theorem~\ref{faithful} we also obtain that
${\rm Hom}(\delta(I), H_1) = R/I_1(\varphi)$. In addition, if
$I_1(\varphi)$ is Cohen-Macaulay of codimension $g$ then by the
theorem of Hartshorne-Ogus we have that $\delta(I)$ $($which is
$S_2)$ is Cohen-Macaulay and therefore ${\rm depth}\,I/I^2 \geq
d-g-2$.}
\end{Remark}

\medskip

Unfortunately, there is not much hope to stretch the proof of
Theorem~\ref{faithful} as the following example shows.

\begin{Example}
{\rm Let $R$ be the localized polynomial ring $k[x,y]_{(x,y)}$.
The ideal $I=(x^5-y^5, x^4y, xy^4)$ is such that $I^2 \colon I =
(I, x^3y^3)$. In this case $I_1(\varphi) = (x, y)^2$ so that
$\mu(I) = \mu(I_1(\varphi)) = 3$. }
\end{Example}

\smallskip

\section{More Questions}
\medskip\
We end by considering some other closely related questions which
came up during the course of this investigation. We let $I$ be an
${\mathfrak m}$-primary ideal of the local ring $R$ minimally generated by
$n$ elements, and let $J_i$ be the
annihilator of the $i$th Koszul cohomology of $I$ with respect to
a minimal generating set of $I$.

\begin{quote}
\emph{Set $d$ equal to the dimension of $R$. Is $J_1\cdot J_2\cdots J_{n-d}$ contained
in $\overline{I^{n-d}}$?}
\end{quote}

Notice that the Koszul homology of $I$ vanishes for values larger than $n-d$,
so that the product above represents all the interesting annihilators of the
Koszul homology of $I$. Furthermore, a postive answer to this question gives a
positive answer to our main question. This follows since each $J_i$ contains $I$.
Along any discrete valuation $v$, this means that $v(I)\geq v(J_i)$ for all
$i$. A positive answer to the question above implies that
$$\sum_{i =1}^{n-d} v(J_i)\geq (n-d)v(I)\geq \sum_{i = 1}^{n-d} v(J_i).$$
It would follow that $v(J_i) = v(I)$ for all $i$, implying that $J_i\subseteq \overline{I}$
for all $1\leq i\leq n-d$. Conversely, if $J_i\subseteq \overline{I}$
for all $1\leq i\leq n-d$, then clearly $J_1\cdot J_2\cdots J_{n-d}$ is contained
in $\overline{I^{n-d}}$, so the above question is equivalent to saying that
$J_i\subseteq \overline{I}$
for all $1\leq i\leq n-d$. This form of the question suggests using homotopies
to compare the Koszul complex of a set of generators of $I$ with the free resolution
of $I$. However, we have not been able to use this idea to settle the question.

Another question which arose during our work is the following:

\begin{quote}
\emph{Let $n$ be the number of minimal generators of an ${\mathfrak m}$-primary ideal
$I$ in a Cohen-Macaulay local ring $R$ with infinite residue field, and
let $d$ be the dimension of the ring.
For every $j$, $d\leq j\leq n-1$, choose $j$ general minimal generators of $I$,
and let $J_j$ be the ideal they generate. Let $H_{n-j}$ denote the $(n-j)$th Koszul homology
of a minimal set of generators of $I$. Is
$${\rm Ann}(H_{n-j})\subseteq J_j:(J_j:I)\subseteq \overline{I}?$$}
\end{quote}

We have positive answers to this question for the two extremes:
$j = d$ and $j = n-1$, in the latter case assuming $I$ is integrally closed.

\begin{center}
{\Large\bf Acknowledgments}
\end{center}
C. Huneke and W.V. Vasconcelos gratefully acknowledge partial
support from the National Science Foundation.

\end{document}